\documentclass[aop,citesort,MSNbibl,dvips]{arximspdf}
\usepackage{graphicx}

\doi{10.1214/10-AOP596}
\volume{39}
\issue{5}
\pubyear{2011}
\firstpage{1983}
\lastpage{2017}

\makeatletter

\newtheorem{theorem}{Theorem}
\newtheorem{lemma}{Lemma}
\newtheorem{cor}{Corollary}
\newtheorem{prop}{Proposition}

\newcommand{\xrightarrow}[1]{\stackrel{#1}{\rightarrow}}

\newcommand{\eps}{\varepsilon}
\newcommand{\bz}{\mathbf0}
\newcommand{\Z}{\mathbb Z}
\newcommand{\R}{\mathbb R}
\newcommand{\C}{\mathbb C}
\newcommand{\Tr}{\operatorname{Tr}}
\newcommand{\SL}{\mathrm{SL}_2\C}
\newcommand{\SU}{\mathrm{SU}_2\C}
\newcommand{\cof}{\operatorname{cof}}
\newcommand{\E}{\mathbb E}
\newcommand{\Ch}{\operatorname{Ch}}
\newcommand{\G}{\mathcal G}
\newcommand{\Qdet}{\operatorname{Qdet}}
\newcommand{\I}{\mathrm{I}}
\newcommand{\trh}{\widehat{\operatorname{tr}}}

\makeatother

\begin{document}
\begin{frontmatter}

\title{Spanning forests and the vector bundle Laplacian\thanksref{T1}}
\runtitle{Spanning forests and the vector bundle Laplacian}

\thankstext{T1}{Supported by the NSF.}

\begin{aug}
\author[A]{\fnms{Richard} \snm{Kenyon}\corref{}\ead[label=e1]{rkenyon@math.brown.edu}}
\runauthor{R. Kenyon}
\affiliation{Brown University}
\dedicated{Dedicated to the memory of Oded Schramm}
\address[A]{Mathematics Department\\
Brown University\\
151 Thayer st.\\
Providence, Rhode Island 02912\\
USA\\
\printead{e1}} 
\end{aug}

\received{\smonth{1} \syear{2010}}
\revised{\smonth{4} \syear{2010}}

%
\begin{abstract}
The classical matrix-tree theorem relates the determinant of the
combinatorial Laplacian on a graph to the number of spanning trees. We
generalize this result to Laplacians on one- and two-dimensional vector
bundles, giving a combinatorial interpretation of their determinants in
terms of so-called cycle rooted spanning forests (CRSFs). We construct
natural measures on CRSFs for which the edges form a determinantal
process.

This theory gives a natural generalization of the spanning tree process
adapted to graphs embedded on surfaces. We give a number of other
applications, for example, we compute the probability that a
loop-erased random walk on a planar graph between two vertices on the
outer boundary passes left of two given faces. This probability cannot
be computed using the standard Laplacian alone.
\end{abstract}

%
\begin{keyword}[class=AMS]
\kwd{82B20}.
\end{keyword}
\begin{keyword}
\kwd{Discrete Laplacian}
\kwd{quaternion}
\kwd{spanning tree}
\kwd{resistor network}.
\end{keyword}

\end{frontmatter}

\section{Introduction}
The classical matrix-tree theorem, which is usually attributed to
either Kirchhoff \cite{Kirchhoff} or Brooks et al.
\cite{BSST}, states that
the product of the nonzero eigenvalues
of the combinatorial Laplacian on a connected graph
is equal to the number of rooted spanning trees of that graph.
This theorem has extensions to graphs with
weighted edges (resistor networks) and more generally
to Markov chains. Forman \cite{Forman}
extended the theorem in a more interesting direction,
to the setting of a line bundle on a graph.
In this paper, we reprove his result and extend
the theorem further to two-dimensional vector bundles on graphs
(with $\SL$ connection).
Our main theorem
gives a combinatorial interpretation of the
determinant of the Laplacian as a sum over
cycle-rooted spanning forests (CRSFs).
These are simply collections
of edges each of whose connected components has as many vertices
as edges (and therefore, each component is a tree plus an edge: a~unicycle).
Here the weight of a configuration is the product of $2$ minus the trace
of the holonomy around each of the cycles.

By varying the connection on the vector bundle, one constructs in this
way many
natural measures on CRSFs. If the monodromy of the underlying connection
is unitary, then these measures are determinantal processes for the edges
(or $q$-determinantal in the case of an $\SL$-connection).

We give here a number of applications of these results.
Further applications can be found in the papers \cite{GK,K2}
and \cite{KW2}.

\subsection{Spanning trees}
One of the main applications is to the study of
spanning trees, in particular spanning trees
on planar graphs or graphs on surfaces.
By a judicious choice of vector bundle connection,
the natural probability measure on CRSFs
can be made to model the uniform spanning tree measure
conditioned on having certain ``boundary connections.'' As simple examples
we compute the probability that, on a finite planar graph,
the branch of the uniform spanning tree connecting
two boundary points $z_1$ and $z_2$
passes left of a given face or a given two faces. These probabilities
are given in terms of the Green's function for the standard Laplacian
on the graph.

\subsection{\texorpdfstring{$\theta$-functions}{theta-functions}}
Suppose $\G$ is a graph embedded on an annulus. On such a graph,
a \textit{$\theta$-function} is a multi-valued harmonic
function whose ``analytic continuation'' around the annulus
is a constant times the original function. In other terminology
(defined below),
it is a harmonic
section of a flat line bundle. The constant is called the
\textit{multiplier} of the $\theta$ function. We prove that
if $k$ is the largest integer such that one can simultaneously
embed $k$ pairwise disjoint cycles in $\G$, each winding once around
the annulus,
then $\G$ has at most $2k-1$ $\theta$-functions. These $\theta
$-functions have
distinct positive real multipliers; these multipliers are closed
under inverses (if $\lambda$ is a multiplier then so is $\lambda^{-1}$).
For generic conductances, $\G$ will have exactly $2k-1$ $\theta$-functions;
we conjecture that this is always the case.
The multipliers have a~probabilistic meaning for CRSFs: see below.

%

\subsection{Graphs on surfaces}
For a graph embedded
on a surface (in such a~way that complementary components are contractible
or peripheral annuli),
a natural probability model is the uniform random CRSF
whose cycles are topologically nontrivial (not null-homotopic). Such CRSFs
are called \textit{incompressible}. The cycles in a CRSF give a finite
lamination of the surface,
that is, a finite collection of disjoint simple closed curves.
The Laplacian determinant on a flat line bundle on the surface counts
CRSFs weighted by a~function of the monodromy of the connection. By considering
the Laplacian determinant as a function on the representation variety
(consisting of flat connections modulo conjugacy), one can extract the terms
for each possible topological type of finite lamination. In particular,
one can study
the uniform measure on incompressible CRSFs. For example, the
number of components for a uniform incompressible CRSF on a
graph on an annulus is distributed as a sum of a finite number of
independent Bernoulli
random variables, with biases given by $\frac{\lambda}{\lambda+1}$,
where $\lambda$ runs over the $\theta$-function multipliers defined above.
See Corollary \ref{Bernoullicor} below.

Given a Riemann surface one can take a sequence of finer and finer graphs
adapted to the metric so that the potential theory on the graph and on the
Riemann surface agree in the limit.
In this case, one can show \cite{K2}
that the uniform incompressible CRSF has a scaling limit, whose distribution
only depends on the conformal structure of the underlying surface.
This is similar to the theorem of Lawler, Schramm and Werner on conformal
invariance of the uniform spanning tree \cite{LSW}.
We compute here for the annulus and the square torus the
exact distribution of the number and homology class of the cycles
of an incompressible CRSF in the scaling limit.

\subsection{Monotone lattice paths}
The papers \cite{HK} and \cite{KW}
studied CRSFs on a~north/east-directed $m\times n$ grid on a torus.
The distribution of the number of cycles was shown to have a highly nontrivial
structure as a function of~$m,n$. Here the line bundle Laplacian gives
quantitative
information about the model which was not available in
\cite{HK}. In particular, we obtain an exact expression for the generating
function for the number and homology type of the cycles.
For the $n\times n$ torus, for example, we show that the number of
cycles tends to a
Gaussian as $n\to\infty$ with expectation $
\sqrt{\frac{n}{4\pi}}$
and variance
$\sqrt{\frac{n}{4\pi}}(1-\frac1{\sqrt{2}})$.

\section{Background}

The uniform probability
measure on spanning trees on a~graph, called the UST measure,
has for the past $20$ years been a remarkably successful
and rich area of study in probability theory.
Pemantle \cite{Pemantle}
showed that the unique path between two points in
a uniform spanning tree has the same distribution as the loop-erased
random walk between those two points.\setcounter{footnote}{1}\footnote{The
loop-erased random walk (LERW)
between $a$ and $b$ is defined as follows:
draw the trace a simple random walk from $a$ stopped when it reaches $b$,
and then erase from the trace all loops in chronological order.
What remains is a simple path from $a$ to $b$.}
Wilson \cite{Wilson}
extended this to give a simple method of sampling
a uniform spanning tree in any graph.
Burton and Pemantle \cite{BP}
proved that the edges of the uniform spanning tree form
a determinantal process (see definition below). This allows
computation of multi-edge probabilities in terms of the Green's function.

The UST on $\Z^2$ received particular attention due to
the conformal invariance properties of its scaling limit.
Pemantle showed that almost surely the UST in $\Z^2$ has one component.
In \cite{K54}, we
proved that the expected length of the LERW in~$\Z^2$,
and therefore a branch of the UST
of diameter~$n$ was of order $n^{5/4}$,
a result predicted earlier by conformal
field theory \cite{Maj54}. In \cite{KenyonCI},
we showed that the scaling limit of the ``winding field''
(describing how the branches of the UST wind around faces)
was a Gaussian free field.
Further conformally invariant properties
were proved in \cite{Klongrange}.
In \cite{LSW}, Lawler, Schramm and Werner proved that the $\Z^2$-LERW
converges to $\operatorname{SLE}_2$,
and the peano curve winding around the $\Z^2$-UST converges
to $\operatorname{SLE}_8$.

As a result of these works, we have a decent understanding
of the scaling limit of the UST on $\Z^2$. However, some
important questions
remain. For example,
how are different points in the UST connected:
given a set of points in the plane, what is the topology of
their tree convex hull, that is, the union of the branches of the UST
connecting them in pairs?
Can one compute various connection probabilities,
for example, the probability that the LERW from $(0,0)$
passes through the points $v_1,v_2,\ldots, v_n$ in order?
What is the distribution of the ``bush size,'' the finite component
of the tree obtained by removing a single random edge? While many of these
can be answered in principle using SLE techniques, in practice one must solve
a~hard PDE.

Many of these questions can be answered
with the bundle Laplacian. Some of these are illustrated below.
We will discuss how these
results can be used to compute various connection
probabilities for the UST in \cite{KW2},
and discuss the connection with integrable systems in \cite{GK}.

\section{Vector bundles on graphs}

\subsection{Definitions}
Let $\G$ be a finite graph.
Given a fixed vector space $V$,
a \textit{$V$-bundle}, or simply a \textit{vector bundle
on} $\G$ is the choice of a vector space $V_v$ isomorphic to $V$ for
every vertex $v$ of $\G$. A vector bundle can be identified with the
vector space
$V_{\G}:=\bigoplus_{v}V_v \cong V^{|\G|}$.
A \textit{section} of a vector bundle is an element of $V_{\G}$.

A \textit{connection} $\Phi$ on a $V$-bundle
is the choice for each edge $e=vv'$ of $\G$ of an isomorphism $\phi_{vv'}$
between the corresponding vector spaces $\phi_{vv'}\dvtx V_v\to V_{v'}$,
with the property that $\phi_{vv'}=\phi_{v'v}^{-1}$. This isomorphism
is called the \textit{parallel transport} of vectors in $V_v$ to vectors
in $V_{v'}$. Two connections $\Phi,\Phi'$ are said to be
\textit{gauge equivalent} if there is for each vertex an isomorphism
$\psi_v\dvtx V_v\to V_v$ such that the diagram\vspace*{6pt}

\begin{center}
\begin{tabular}{c}

\includegraphics{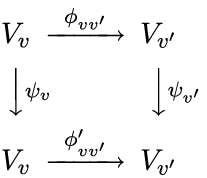}
\vspace*{6pt}
\end{tabular}
\end{center}

\noindent commutes. In other words $\Phi'$ is just a base change of $\Phi$.
Given an oriented cycle $\gamma$ in $\G$ starting at $v$, the \textit{monodromy} of the connection is
the\vadjust{\eject}
element of $\operatorname{End}(V_v)$ which is the product of the
parallel transports
around $\gamma$.
Monodromies starting at different vertices on $\gamma$ are conjugate,
as are monodromies of gauge-equivalent connections.

A \textit{line bundle} is a $V$-bundle where $V\cong\C$,
the one-dimensional complex vector space. In this case, if we choose a
basis for each $\C$
then the parallel transport is just multiplication by an element of $\C
^*=\C\setminus\{0\}$.
The monodromy of a cycle is in $\C^*$ and does not depend
on the starting vertex of the cycle (or gauge).

\subsection{The Laplacian}

The Laplacian $\Delta$ on a $V$-bundle with connection
$\Phi$ is the linear operator
$\Delta\dvtx V_{\G}\to V_{\G}$ defined by
\[
\Delta f(v)=\sum_{v'\sim v}f(v)-\phi_{v'v}f(v'),
\]
where the sum is over neighbors $v'$ of $v$.\vadjust{\goodbreak}

If we assign to each edge a positive real weight (a conductance)
$c_{vv'}=c_{v'v}$, the weighted Laplacian is defined by
\[
\Delta f(v)=\sum_{v'}c_{vv'}\bigl(f(v)-\phi_{v'v}f(v')\bigr).
\]

Note that if the vector bundle is \textit{trivial}, in the sense that
$\phi_{vv'}$ is the identity
for all edges, this is the classical notion of graph Laplacian (or more
precisely, the direct sum of $\dim V$ copies of the Laplacian).

Here is an example.
Let $\G=K_3$ with vertices $\{v_1,v_2,v_3\}$.
Let $\Phi$ be the line bundle connection
with $\phi_{v_iv_j}=z_{ij}\in\C^*$.
Then in the natural basis, $\Delta$ has matrix
%
\begin{equation}\label{3vexample}
\Delta=\pmatrix{2&-z_{12}&-z_{13}\cr-z_{12}^{-1}&2&
-z_{23}\cr-z_{13}^{-1}&-z_{23}^{-1}&2}.
\end{equation}

\subsection{Edge bundle}

One can extend the definition of a vector bundle to the edges of
$\G$. In this case, there is a vector space $V_e\cong V$ for each edge
$e$ as well as each vertex.
One defines connection isomorphisms $\phi_{ve}=\phi_{ev}^{-1}$ for a
vertex $v$
and edge $e$ containing that vertex, in such a way that if $e=vv'$ then
$\phi_{vv'}=\phi_{ev'}\circ\phi_{ve}$, where $\phi_{vv'}$ is
the connection on the vertex bundle.

The vertex/edge bundle can be identified with $V^{|\G|+|E|}=V_{\G}\oplus V_E$,
whe\-re~$V_E$ is the direct sum of the edge vector spaces.

A \textit{$1$-form} (or cochain) is a function on oriented edges which is
antisymmetric
under changing orientation. If we fix an orientation for each edge,
a~$1$-form is a section of the edge bundle, that is, an element of $V^{|E|}$.
We denote by $\Lambda^1(\G,\Phi)$ the space of $1$-forms
and $\Lambda^0(\G,\Phi)$ the space of $0$-forms, that is, sections
of the vertex bundle.

We define a map $d\dvtx\Lambda^0(\G,\Phi)\to\Lambda^1(\G,\Phi)$
by $df(e)=\phi_{ye}f(y)-\phi_{xe}f(x)$ where $e=xy$ is an oriented
edge from vertex
$x$ to vertex $y$.
We also define
an operator $d^*\dvtx\Lambda^1\to\Lambda^0$ as follows:
\[
d^*\omega(v) = \sum_{e=v'v}\phi_{ev}\omega(e),
\]
where the sum is
over edges containing $v$ and oriented toward $v$. Despite the notation,
this operator $d^*$ is not a standard
adjoint of $d$ unless $\phi_{ev}$ and $\phi_{ve}$ are adjoints
themselves, that is, if parallel transports are unitary operators (see below).

The Laplacian $\Delta$ on $\Lambda^0$ can then be defined as
the operator $\Delta=d^*d$:
\begin{eqnarray*}
d^*\,df(v)&=& \sum_{e=v'v}\phi_{ev}\,df(e)\\
&=&\sum_{e=v'v} \phi_{ev}\bigl(\phi_{ve}f(v)-
\phi_{v'e}f(v')\bigr)\\
&=&\sum_{v'}f(v)-\phi_{v'v}f(v')\\
&=&\Delta f(v).
\end{eqnarray*}

We can see from the example (\ref{3vexample}) above on $K_3$
that $\Delta$ is not necessarily self-adjoint.
However, if $\phi_{vv'}$ is unitary: $\phi_{vv'}^{-1}=\phi_{v'v}^*$
then $d^*$ will be the adjoint of $d$
for the standard Hermitian inner products
on $V^{|\G|}$ and $V^{|E|}$, and so in this case
$\Delta$ is a Hermitian, positive semidefinite
operator.
In particular on a line bundle if $|\phi_{vv'}|=1$ for all
edges $e=vv'$, then
$\Delta$ is Hermitian and positive semidefinite.

\section{Spanning forests associated to a line bundle}

\subsection{Cycle-rooted spanning forests (CRSFs)}

\begin{figure}
\includegraphics{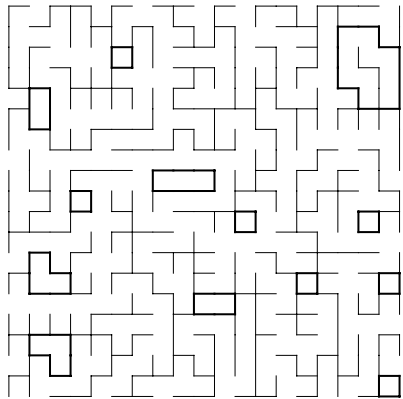}
\caption{A CRSF on a square grid.}
\label{CRSF}
\vspace*{-3pt}
\end{figure}

For a line bundle $(\G,\Phi)$ on a finite graph $\G$, we have a combinatorial
interpretation of the determinant of the Laplacian in terms of
cycle-rooted spanning forests.
A \textit{cycle-rooted spanning forest} (CRSF) is a subset $S$ of the
edges of $\G$,
spanning all vertices (in the sense that every vertex is the endpoint
of some edge)
and with the
property that each connected component of $S$
has as many vertices as edges (and so has a unique cycle). See Figure
\ref{CRSF}.

An \textit{oriented CRSF} is a CRSF in which we orient each edge
in such a way that cycles are oriented coherently
and we orient the edges not in a cycle
toward the cycle. An oriented CRSF is the same as a
\textit{nonzero vector field} on $\G$, which is the choice of a single
outgoing edge from each vertex.

A \textit{cycle-rooted tree} (CRT), also called unicycle,
is a component of a CRSF (oriented or not).

\subsection{The matrix-CRSF theorem}

The following theorem is due to Forman~\cite{Forman}.
\begin{theorem}[\cite{Forman}]\label{mainline} For a line bundle on a
connected finite graph,
\[
\det\Delta=\sum_{\mathrm{CRSFs}}\prod_{\mathrm{cycles}}(2-w-1/w),
\]
where the sum is over all unoriented CRSFs $C$,
the product is over the cycles of $C$, and $w, 1/w$ are the
monodromies of the two orientations of the cycle.
\end{theorem}

Recall that for a line bundle
the monodromy of an oriented cycle does not depend on the
starting vertex.
Note that we could as well have written
$\det\Delta=\sum_{\mathrm{OCRSFs}}\prod_{\mathrm{cycles}}(1-w)$
where now the sum is over oriented CRSFs (OCRSFs).

Forman uses an explicit
expansion of the determinant as a sum over the symmetric group,
and a careful rearrangement of the terms. We give a different proof using
Cauchy--Binet formula which allows us to bypass this step.\looseness=-1
\begin{pf*}{Proof of Theorem \ref{mainline}}
We use the Cauchy--Binet formula,
$\det\Delta=\det d^*d=\sum_B \det(B^*B)$, where $B$ is a maximal
minor of $d$ (i.e., if $d$ is $n\times m$ then $B$ runs
over all $n\times n$
submatrices of $d$ obtained by choosing $n$ columns).

Columns of $B$ index vertices and rows of $B$ index edges of $\G$.
Since $B$ is square, it
corresponds to a selection of $n$ edges (where $n$ is the number of
vertices of $\G$).

If the set of edges in $B$ contains a component $C$ with no cycles, say
of~si\-ze $|C|=k$, this component necessarily has $k-1$ edges.
So $B$ maps a $k$-dimen\-sional subspace of the vertex space
to a $k-1$-dimensional
subspace of the edge space and so $\det B=0$. Therefore to have a
nonzero contribution
to the sum every component of $B$\vadjust{\eject}
has at least one cycle, and, since the total number of vertices equals
the total
number of edges, every component must have exactly one cycle.
We have proved that the sum is over CRSFs.

If $B$ has only one component which is a CRT, orient its
edges so that all edges are
oriented toward the unique cycle, and along the cycle choose one of
the two coherent orientations.
Write
\[
\det B = \sum_{\sigma\in S_n}\operatorname{sgn}(\sigma)B_{1\sigma(1)}
\cdots B_{n\sigma(n)}.
\]
The only nonzero
terms in this expansion are those in which
vertex $i$ is adjacent to edge $\sigma(i)$, that is,
the two terms which are consistent with one of these
two orientations.
Thus
%
\begin{eqnarray}\label{line2}
\det B&=&\operatorname{sgn}(\sigma_1)\biggl(\prod_{j\in\mathrm
{bushes}} -\phi_{v_je_j}\prod_{i\in\mathrm{cycle}}-\phi
_{v_ie_i}\biggr)+\operatorname{sgn}(\sigma_2)\nonumber\\[-8pt]\\[-8pt]
&&{}\times\biggl(\prod_{j\in
\mathrm{bushes}} -\phi_{v_je_j}
\prod_{i\in\mathrm{cycle}}-\phi_{v_{i+1}e_i}\biggr)\nonumber\\
&=&
\operatorname{sgn}(\sigma_1)(-1)^n
\prod_{j\in\mathrm{bushes}} \phi_{v_je_j}\biggl(\prod_{i\in
\mathrm{cycle}}\phi_{v_ie_i}-
\prod_{i\in\mathrm{cycle}}\phi_{v_{i+1}e_i}\biggr),\nonumber
\end{eqnarray}
where the products over $j$ are over edges not in the cycle
(called bush edges),
and those over $i$ are for the two cyclic
orientations (indices are taken cyclically).
If the cycle has odd length, $\operatorname{sgn}(\sigma
_1)=\operatorname{sgn}(\sigma_2)$ and
the sign change in front of the second product in (\ref{line2})
is due to switching the sign on the odd number of edges of the cycle.
If the cycle has even length, the sign change is due to $\operatorname
{sgn}(\sigma_1)=
-\operatorname{sgn}(\sigma_2)$ (a~cyclic permutation of even length
has signature
$-1$).

We similarly have
\[
\det B^*=\operatorname{sgn}(\sigma_1)(-1)^n
\prod_{j\in\mathrm{bushes}} \phi_{e_jv_j}\biggl(\prod_{i\in\mathrm
{cycle}}\phi_{e_iv_i}-
\prod_{i\in\mathrm{cycle}}\phi_{e_iv_{i+1}}\biggr),
\]
and multiplying these two and using $\phi_{ev}\phi_{ve}=1$, we get
\[
\det B^*B=(A_1-A_2)\biggl(\frac1{A_1}-\frac1{A_2}\biggr) = 2-\frac
{A_1}{A_2}-\frac{A_2}{A_1},
\]
where $\frac{A_1}{A_2}=\prod_{i\in\mathrm{cycle}}\phi
_{v_iv_{i+1}}$\vspace*{2pt} which is the monodromy of
the cycle.

If each component of $B$ is a CRT, then $\det B^*B$ is the product over
components
of the above contributions.
\end{pf*}

\section{Measures on CRSFs}

In the case of a line bundle $\Phi$ with unitary connection, Theorem
\ref{mainline} provides a definition of a natural probability
measure~$\mu_{\Phi}$ on CRSFs.
A CRSF has probability equal to
$\frac1Z\prod(2-w-\frac1w)$ where the product is over its cycles.
The constant of proportionality $Z$, or partition function,\vadjust{\eject}
is then just $\det\Delta$. Note that as long as the (unitary) bundle
is not
trivializable (gauge equivalent to the trivial bundle) then some
cycle has nontrivial monodromy and so $\det\Delta>0$.

The goal of this section is to show that the measure is determinantal
for the edges.

\subsection{Harmonic sections}

The Laplacian is the composition:
\[
\Lambda^0\xrightarrow{d}\Lambda^1\xrightarrow{d^*}\Lambda^0.
\]
If $\Delta$ has full rank, then $\Lambda^1(\G,\Phi)$ splits as
\[
\Lambda^1=\operatorname{Im}(d)\oplus\operatorname{Ker}(d^*).
\]

A section $f\in\Lambda^0$ is \textit{harmonic} if $\Delta f=0$.
If $\Delta$ does not have full rank, there are nontrivial harmonic sections
and
the above is no longer a direct sum.
A $1$-form $\omega\in\Lambda^1$ is harmonic if it is in
the intersection $\operatorname{Im}(d)\cap\operatorname{Ker}(d^*)$.
Such a form is both
\textit{exact} ($\omega=df$ for some $f$) and \textit{co-closed} ($d^*\omega=0$).
\begin{prop}\label{Pdef}
When $\Delta$ has full rank,
the projection along $\operatorname{Ker}(d^*)$ from $\Lambda^1$ to the space of exact
$1$-forms $\operatorname{Im}(d)$
is given by the operator
$P=dgd^*$ where $g=\Delta^{-1}$ is the Green's function for $\Delta$
on $\Lambda^0$.
\end{prop}
\begin{pf} Note that $P=dgd^*$ is zero on co-closed forms
and the identity on exact forms: $P(df) = dgd^*df=dg\Delta f=df$.
\end{pf}

This projection operator $P$ plays a role below.

\subsection{Determinantal processes}

Let $\mu$ be a probability measure on $\Omega=\{0,1\}^n$.
For a set of indices $S\subset\{1,2,\ldots,n\}$,
define the set
\[
E_S=\{(x_1,\ldots,x_n)\in\Omega\mid x_s=1\ \forall s\in S\}.
\]
We say $\mu$ is a \textit{determinantal} probability measure (see, e.g.,
\cite{Peresetal})
if there is an $n\times n$ matrix $K$ with the property that
for any set of indices $S\subset\{1,2,\ldots,n\}$,
we have
\[
\mu(E_S)=\det[(K_{i,j})_{i,j\in S}].
\]
That is, principal minors of $K$ determine the probability of events
of type~$E_S$.

The matrix $K$ is called the \textit{kernel} of the measure.
A simple example of a~determinantal process
is the product measure, in which case $K$ is a~diagonal matrix
with diagonal entries in $[0,1]$.
Another well-known example of a determinantal process is
the uniform spanning tree measure on a finite graph~\cite{BP}.
Here $\{1,2,\ldots,n\}$ index the edges of a connected graph.
The kernel~$K$ is the so-called \textit{transfer-current matrix},
defined by $K=dGd^*$, where~$G$ is the Green's function,
see \cite{BP}.

A simple inclusion-exclusion argument shows that
individual point probabilities are also given by determinants:
\begin{prop}\label{pointprobs}
Let $X=(x_1,x_2,\ldots,x_n)\in\Omega$. Then
\[
\mu(X)= (-1)^{n-|X|} \det\bigl(\operatorname{diag}(X) - K\bigr),
\]
where $\operatorname{diag}(X)$ is the diagonal matrix whose diagonal entries
are $1-x_i$, and $|X|=\sum_{i=1}^n x_i$.
\end{prop}

\subsection{A determinantal process on CRSFs}

We will need the following well-known
lemma.
\begin{lemma}\label{Schur} If $\bigl({{A\atop C} \enskip{B\atop D}}\bigr)$ is a
block matrix
with $A,D$ square and $D$ invertible then
\[
\det\pmatrix{A&B\cr C&D}=\det D\det(A-BD^{-1}C).
\]
\end{lemma}
\begin{pf} This follows from
\[
\pmatrix{A&B\cr C&D}=\pmatrix{A-BD^{-1}C&BD^{-1}
\cr0&I}\pmatrix{I&0\cr C&D}.
\]
\upqed\end{pf}

Recall the definition of the projection $P$ from $1$-forms to exact $1$-forms
(Proposition \ref{Pdef}).
\begin{theorem}
For a line bundle $\Phi$ with unitary connection, with respect
to the measure $\mu_\Phi$ the edges of the CRSF form
a determinantal process with kernel $P=dgd^*$.
\end{theorem}
\begin{pf}
Let $\G$ be a graph with $n$ vertices and $m$ edges.
Let $\{e_1,\ldots,e_n\}\subset E$ be the edges of a CRSF $\gamma$.
Write the matrix for $d$ and $d^*$ so that the first~$n$ edges are
$e_1,\ldots,e_n$.
Then $d={d_1\choose d_2}$ where $d_1$ is $n\times n$
and $d_2$ is $(m-n)\times n$.

From Proposition \ref{pointprobs}, the probability of $\gamma$ is
$\Pr(\gamma) = (-1)^{m-n}\det(X-dgd^*)$ which by Lemma \ref{Schur}
with $A=X, B=d, D=\Delta$ and $C=d^*$ can be written
\[
\Pr(\gamma) = (-1)^{m-n}\frac{\det\left({{0\atop 0}\atop d_1^*}\enskip{{0\atop I_{m-n}}\atop d_2^*}\enskip
{{d_1\atop d_2}\atop \Delta}
\right)}{\det\Delta} = \frac{\det d_1^*d_1}{\det\Delta}.
\]
By Theorem \ref{mainline} this is exactly the probability of $\gamma$.
\end{pf}

Among other things, this theorem along with Proposition \ref{pointprobs}
allows us to do exact sampling from the measure $\mu_\Phi$, as follows
(see \cite{Peresetal}).
Pick an edge~$e_1$;
it is present with probability $P(e_1,e_1)$.
Take another edge~$e_2$; if $e_1$ is present,~$e_2$ will be present
with probability
\[
\Pr(e_2|e_1) = \frac{1}{P(e_1,e_1)}\det\pmatrix{
P(e_1,e_1)&P(e_1,e_2)\cr P(e_2,e_1)&P(e_2,e_2)}.
\]
If $e_1$ is not present,
$e_2$ will be present with probability
\[
\Pr(e_2| \neg e_1)= \frac{-1}{1-P(e_1,e_1)}
\det\left(\pmatrix{1&0\cr0&0}-\pmatrix{
P(e_1,e_1)&P(e_1,e_2)\cr P(e_2,e_1)&P(e_2,e_2)}\right)
\]
and so on.

\subsection{Small monodromies}
An important limit, or actually set of limits of these determinantal
processes, is when
the monodromies tend to $1$. Even though $\det\Delta$
tends to zero, the probability measures $\mu_\Phi$
may converge.
We get determinantal
processes supported on CRSFs on the graph with trivial line bundle.
\begin{theorem}
Choose an orientation for each edge and
suppose the parallel transport is $\phi_{e_j}=e^{itc_j}$ for the $j$th
edge in the direction of its orientation.
Fix the $c_j$ and let $t\to0$.
In the limit $t\to0$ the determinantal process $\mu_{\Phi_t}$ tends
to a determinantal process $\mu_0=\mu_0(\{c_j\})$ supported
on CRSFs with a~single component,
that is, CRTs.
The probability of a CRT is proportional to $(\sum c_j)^2$,
where the sum is over oriented edges in the unique cycle.
\end{theorem}
\begin{pf}
The monodromy of a loop is $w=e^{it\sum c_j}$. We have $2-w-\frac
1w=2-2\cos(t\sum c_j)=t^2(\sum c_j)^2+O(t^3)$. In the limit $t\to0$,
the partition function is $Z=O(t^2)$ so only CRSFs with one component
contribute.
\end{pf}

\section{\texorpdfstring{An application to $\G\times\Z_n$}{An application to G times Z_n}}

Let $\G$ be a finite graph with $m$ vertices
and $\Z_n$ the $n$-cycle. Let $H_n$ be the product graph whose
vertices are
$\G\times\Z_n$ and edges connect $(x,j)$ to $(x,j+1\bmod n)$ and
$(x,j)$ to $(x',j)$ when $x,x'$
are neighbors in $\G$.

We compare the uniform spanning tree on $H_n$ and
the uniform cycle-rooted tree whose cycle winds around $\Z_n$.
The \textit{minimum cut set} of a spanning tree is the set of edges not in
the tree
such that when added to the tree make a cycle winding nontrivially
around $\Z_n$.

Fix $z\in\C$ with $|z|=1$.
On the product graph $H_n$ let $\phi_{vv'}=1$ except when
$v=(x,n-1)$ and $v'=(x,0)$ in which case $\phi_{vv'}=z$. Then
cycles have trivial monodromy unless they wind nontrivially around the
$\Z_n$ direction.
So a CRSF has nonzero probability only if all of
its cycles wind nontrivially
around $\Z_n$ (possibly many times).

Because of the product structure the eigenvectors of the Laplacian on
$H_n$ are
products of eigenvectors of the standard Laplacian on $\G$ and the
line bundle Laplacian on $\Z_n$. Thus, the eigenvalues
of the Laplacian on $H_n$
are sums of eigenvalues of $\Delta_{\G}$ (the standard Laplacian)
and
$\Delta_{\Z_n}$ (the line bundle Laplacian).
The eigenvalues of $\Delta_{\Z_n}$
are $2-\zeta-\zeta^{-1}$ where $\zeta$ ranges over the roots
of $\zeta^n=z$ (the corresponding eigenvectors are exponential
functions).\vadjust{\eject} The Laplacian determinant of $H_n$ is then
%
\begin{eqnarray}\label{Ncrsfs}
\det\Delta&=& \prod_{\zeta^n=z} \det(\Delta_{\G}+2-\zeta-\zeta
^{-1})=\prod_{\zeta^n=z}\prod_\lambda
(\lambda+2-\zeta-\zeta^{-1})
\nonumber\\[-8pt]\\[-8pt]
&=&(2-z-z^{-1})\prod_{\zeta^n=z}\prod_{\lambda\ne0}(\lambda+2-\zeta
-\zeta^{-1}),\nonumber
\end{eqnarray}
where $\lambda$ runs over the eigenvalues of the Laplacian on $\G$.

Compare this to $\det '\Delta_0$,
the product of nonzero eigenvalues
of the standard Laplacian $\Delta_0$ on $H_n$:
%
\begin{eqnarray}\label{nrootedtrees}
\det{'}\Delta_0&=&\prod_{(\zeta,\lambda)\ne(1,0)}\lambda+2-\zeta
-\zeta^{-1}\nonumber\\
&=&\prod_{\zeta^n=1,\zeta\ne1} 2-\zeta-\zeta^{-1}\prod_{\zeta^n=1}
\prod_{\lambda\ne0}\lambda+2-\zeta-\zeta^{-1}\\
&=&n^2\prod_{\zeta^n=1}\prod_{\lambda\ne0}\lambda+2-\zeta-\zeta
^{-1}.\nonumber
\end{eqnarray}
\begin{theorem}
As $n\to\infty$ a $\mu_{\Phi}$-random CRSF has one component
with probability tending to $1$, that is, is a CRT. Its unique cycle
winds around~$\Z_n$ once with probability tending to $1$.
The ratio $R_n$ of the number\vspace*{1pt} of such CRTs and number of
spanning trees satisfies $\lim_{n\to\infty}\frac{R_n}{n}=\frac1{m}$.
The expected\vadjust{\goodbreak} length~$\E(L_n)$ of the cycle in such a random CRT of $H_n$
satisfies
\[
\lim_{n\to\infty} \frac{\E(L_n)}{n\E(S_n)}=\frac{1}{m},
\]
where $\E(S_n)$ is
the expected size of the minimum cut set of a random spanning tree of $H_n$.
\end{theorem}
\begin{pf}
The first two statements can be seen as follows.
Take $n$ large and consider the part of the
configuration in a piece $\G\times[i,i+N]$ for
large $N$. We claim that
the number of spanning tree configurations restricted to this subgraph,
as a function of $N$, has a strictly
larger exponential growth rate than that of the number of spanning
forest configurations in this subgraph
with two or more components (``strands'') each of which intersects both
ends $\G\times\{i\}$ and
$\G\times\{i+N\}$. The growth rate
of spanning trees can
be computed as the leading eigenvalue of nonnegative finite matrix $T$,
the transfer matrix.\footnote{The states of the transfer matrix
are all possible
forests (edge subsets) of $\G\times\{i\}$ occurring in a spanning
tree of $\G\times\Z$
along with a partition of the components in this forest, which
describes which
components of this forest are connected to each other to the left, that
is, in the part
of the tree in $\G\times(-\infty,i-1]$.}
It is easily seen that the transfer matrix $T$ is primitive, that is,
has the property that some power is strictly positive. On the other
hand, the
transfer matrix for spanning forest configurations with two or more strands
is a~matrix which can be obtained from $T$ by setting some of its positive
entries to $0$
(those entries whose components are all connected to the left).
This strictly decreases its leading eigenvalue. This proves the claim
and the first two statements of the theorem.

If we divide (\ref{Ncrsfs}) by $2-z-z^{-1}$ and let $z$ tend to $1$
the expression counts CRSFs with one component, that is, CRTs,
with a weight
$k^2$ if they wind $k$ times around $\Z/n\Z$, because
\[
\lim_{z\to1}\frac{2-z^k-z^{-k}}{2-z-z^{-1}}=k^2.
\]
As $n$ gets large CRTs which wind more than once around have
probability tending to zero.
Line (\ref{nrootedtrees}) is the number of rooted trees; dividing by
$nm$, the number of locations for the root, gives the number of trees,
which by (\ref{Ncrsfs}) and the above remarks is $\frac{n}{m}$ times
the number of CRTs winding once around,
plus errors tending to zero as $n\to\infty$. This proves the third statement.

The following two sets $A$ and $B$ are in bijection:
$A$ is the set of CRTs whose cycle winds
once around $\Z_n$, along with the choice of an edge on this cycle.
$B$ is the set of spanning trees with a choice of
a complementary edge in the minimum cut set.
The bijection consists in adding the edge to the tree.
Therefore, the expected length of the unique cycle in a random CRT
is
\[
\frac{|A|}{|\mathrm{CRTs}|} = \frac{|B|}{|\mathrm{CRTs}|}=\frac{|B|}{|
\mathrm{trees}|}\frac{|\mathrm{trees}|}{|\mathrm{CRTs}|}.
\]
By the above, $\frac{|\mathrm{trees}|}{|\mathrm{CRTs}|}\to\frac{n}{m}$.
\end{pf}

\section{Weighted and/or directed graphs}

\subsection{Edge weights}

Putting weights on the edges is a minor generalization;
let $c\colon E\to\R_{>0}$ be a conductance associated to each edge,
with $c_{vv'}=c_{v'v}$.
We then define
\[
\Delta f(v) = \sum_{v'} c_{vv'}\bigl(f(v)-\phi_{v'v}f(w)\bigr).
\]

In other words, letting $C$ be the diagonal matrix, indexed by the
edges, whose
diagonal entries are the conductances, we have
\[
\Delta=d^*Cd.\vspace*{-6pt}
\]
\begin{theorem}
$\det\Delta=\sum_{\mathrm{CRSFs}} \prod_{e} C(e)\prod(2-w-\frac1w)$
where the first product is over all edges of the configuration and the second
is over cycles, and~$w$ is the monodromy of the cycle as before.
\end{theorem}

This follows directly from the proof of Theorem \ref{mainline} above.\vadjust{\eject}

\subsection{Directed graphs}

We can also have weights which depend on direction,
so that each edge has two weights, one associated to each direction.
This makes $\G$ into a Markov chain, in which the transition
probabilities are proportional to the edge weights.

In this case the Laplacian has the same form
\[
\Delta f(v) = \sum_{v'} c_{vv'}\bigl(f(v)-\phi_{v'v}f(v')\bigr),
\]
except that $c_{vv'}\ne c_{v'v}$ in general.
The operator
$d\dvtx\Lambda^0\to\Lambda^1$ can be defined by
\[
df(e)= \phi_{ve}f(v)-\phi_{v'e}f(v')
\]
as in the unweighted case,
where $e=v'v$,
and $d^*$ by
\[
d^*(\omega)(v) = \sum_{e=v'v}c_{vv'}\phi_{ev}\omega(e),
\]
where the sum is over edges $e=v'v$ containing $v$.

In this case, each oriented CRSF has a different weight. Forman's
theorem in this
case is
the following.
\begin{theorem}[\cite{Forman}]\label{directed}
\[
\det\Delta=\sum_{\mathrm{OCRSFs}}\prod_{e\in\mathrm
{bushes}}c(e)\prod_{\mathrm{cycles}\ \gamma}C(\gamma)\bigl(1-w(\gamma)\bigr),
\]
where the sum is over oriented CRSFs, the first product is
over the edges in the bushes (i.e., not in the cycles), oriented
toward the cycle, and the second product is over
oriented cycles, $C(\gamma)$ is the product of semiconductances
along $\gamma$ and $w$ is the monodromy
of $\gamma$.
\end{theorem}

An application is in Section \ref{monotone} below.

\section{Laplacian with Dirichlet boundary}

Let $\G$ be a graph with line bundle and connection
$\Phi$ and $S$ a subset of its vertices,
which play the role of boundary vertices. We can define a
Laplacian $\Delta_{\G,S}$ with Dirichlet boundary conditions at $S$
as follows:
for $f\in V^{\G\setminus S}$ and $v\in\G\setminus S$,
\[
\Delta_{\G,S} f(v) =(\deg v)f(v)-\sum_{v'\in\G\setminus S, v'\sim
v}\phi_{v'v}f(v').
\]
This is of course just the Laplacian $\Delta_\G$ restricted to the subspace
$V^{\G\setminus S}$ of functions which are zero on $S$
and projected back to this subspace.
As a~matrix, it is just a submatrix of $\Delta_\G$.

An \textit{essential CRSF} of $(\G,S)$
is a set of edges in which each component is either a
tree containing a unique vertex in $S$ or a cycle-rooted tree
containing no vertices in $S$.
\begin{theorem}\label{theo7}
$\det\Delta_{\G,S}$ is the weighted sum of essential CRSFs,
where each configuration
has weight $\prod_{\mathrm{CRTs}}(2-w-w^{-1})$, where
$w$ is the monodromy of the cycle.
\end{theorem}

This generalizes the matrix-tree theorem\footnote{One
version of the classical matrix-tree theorem is that
the determinant of the matrix obtained by removing
a single row and column $v$ from $\Delta$ is equal to the number of
spanning trees rooted at $v$.}
(when $\Phi$ is the identity and $|S|=1$)
and the matrix-CRSF theorem (when $|S|=0$).
\begin{pf*}{Proof of Theorem \ref{theo7}}
The proof follows the same lines as the proof
of Theorem \ref{mainline}. Now $d$ is an $m\times n$ matrix
with $n=|\G|-|S|$. A maximal minor $B$ of $d$
corresponds to a choice of $n$ edges. Each component
of $B$ is either a CRT, in which case its weight
is computed as before, or a tree connecting some vertices
of $\G\setminus S$ with a single vertex of $S$. In this case,
the determinant of $B$ (we mean, that part of $B$
coming from this component) has a single nonzero term in its
expansion, and the corresponding term in~$B^*$ is its inverse.
So each tree component counts $1$.
\end{pf*}

Again when $|w|=1$ there is a determinantal measure
associated with $\Delta$.

\section{Bundles with $\SL$ connections}

We can extend many of the above results
to the case of a $\C^2$-bundle with $\SL$-connection.

\subsection{$Q$-determinants}

Given $A\in \mathrm{GL}_2(\C)$ define $\tilde A=(\det A)A^{-1}$,
the adjugate of $A$.
That is, if $A=\bigl({{a\atop c} \enskip{b\atop d}}\bigr)$
then $\tilde A=\bigl({{d\atop c} \enskip{b\atop a}}\bigr)$.
Note that
%
\begin{equation}\label{scalar}
A+\tilde A=\operatorname{Tr}(A)I
\end{equation}
is a scalar.

Let $M$ be a matrix with entries in $\mathrm{GL}_2(\C)$.
$M$ is said to be \textit{self-dual} if $M_{ij}=\tilde M_{ji}$.
In particular, diagonal entries must be scalar matrices.

For a matrix $M$ with noncommuting entries,
there are many possible ways one might define its determinant:
however these all involve the expansion
%
\begin{equation}\label{Mdetexp}
\det M=\sum_{\sigma\in S_n}\operatorname{sgn}(\sigma)\bigl[M_{1\sigma(1)}
\cdots M_{n\sigma(n)}\bigr]_{\sigma},
\end{equation}
where
$[M_{1\sigma(1)}
\cdots M_{n\sigma(n)}]_\sigma$ denotes the product of the terms
$M_{i\sigma(i)}$
after they have been rearranged in some specified order depending on
$\sigma$.
Each of these ``order functions'' determines a different possible determinant.

In the current case, there is one very natural condition
to put on these orders. If $\sigma$ is written as a product
of disjoint cycles, then in the corresponding
order the terms in each cycle of $\sigma$
should appear consecutively. Moreover, if $\sigma,\sigma'$
are two permutations with the same cycles, except that some
of the cycles are reversed, then the corresponding order should
have the cycles appear in the same relative order.

The advantage of this for self-dual matrices is that
the product of entries along a cycle is the adjugate of
the product along the reversed cycle; so that the sum of these
is a scalar by (\ref{scalar}). Thus by grouping permutations into
sets with the same cycles up to reversals, one arrives at
a product of scalar matrices. The sum (\ref{Mdetexp}) with
appropriately rearranged
products then yields a~$2\times2$ scalar matrix $q\cdot I$,
and the number $q$ is defined to be the determinant.

The $Q$-determinant of self-dual matrix $M$ \cite{Mehta}
is defined in precisely this way. We define
\[
\Qdet(M) = \sum_{\sigma\in S_n} \operatorname{sgn}(\sigma)
\prod_{\mathrm{cycles}} \frac12\operatorname{tr}(w),
\]
where the sum is over the symmetric group,
each permutation $\sigma$ is written as a product of disjoint cycles,
and $\operatorname{tr}(w)$ is the trace of the product of the matrix entries
in that cycle.
If we group together terms above with the same
cycles---up to the order of traversal of each cycle---then the contribution
from each of these terms
is identical: reversing the orientation of a cycle does not change
its trace.
So we can write
\[
\Qdet(M)=\sum(-1)^{c+n}\prod_{i=1}^c \widehat{\operatorname{tr}}(w_i),
\]
where the sum is over cycle decompositions of the indices, $c$ is the
number of cycles,
and $w_i$ is the monodromy (in one direction or the other)
of each cycle. Here $\widehat{\operatorname{tr}}$ is equal to the
trace for
cycles of length\vspace*{1pt} at least $3$; cycles of length~$1$ or $2$
are their own reversals so we define
$\widehat{\operatorname{tr}}(w)= \frac12\operatorname{tr}(w)$
for these cycles.\looseness=-1

As an example, let $A=a\I, C=c\I$ and $B=\bigl({{b_1\atop b_3} \enskip
{b_2\atop b_4}}\bigr)$. Then
\[
\operatorname{Qdet}\pmatrix{A&B\cr\tilde B&C}=
\trh(A)\trh(C)-\trh(B\tilde B) = ac-(b_1b_4-b_2b_3).
\]

Note that if $M$ is a self-dual $n\times n$ matrix
then $ZM$, considered as a~$2n\times2n$ matrix is antisymmetric,
where $Z$ is the matrix with diagonal blocks
$\bigl({{0\atop-1} \enskip{1\atop0}}\bigr)$ and zeros elsewhere.
The following theorem allows us to compute
$Q$-determinants explicitly.
\begin{theorem}[\cite{Mehta}]
Let $M$ be an $n\times n$ self-dual matrix with entries in $\mathrm{GL}_2(\C)$
and $M'$ the associated $2n\times2n$ matrix, obtained by replacing
each entry with the $2\times2$ block of its entries.
Then $\Qdet(M)=\mathrm{Pf}(Z M')$, the Pfaffian of
the antisymmetric matrix $ZM'$.
\end{theorem}

\subsection{Laplacian determinant}

For a $\C^2$-bundle on a graph $\G$
with $\SL$ connection
$\Delta$ is a self-dual operator (in the above sense).\vadjust{\eject}
\begin{theorem}\label{mainSL2}
For a $\C^2$-bundle on a connected graph
$\G$ with $\SL$ connection,
\[
\Qdet\Delta=\sum_{\mathrm{CRSFs}}\prod_{\mathrm
{cycles}}2-\operatorname{tr}(w),
\]
where $\operatorname{tr}(w)$ is the trace of the monodromy of the cycle.
\end{theorem}
\begin{pf}
The proof generalizes the proof of Theorem \ref{mainline}.
We use as before
the Cauchy--Binet formula which holds for general $Q$-determinants
(Theorem~\ref{CB} below).

Suppose $\G$ has $n$ vertices.
Write $\Qdet\Delta=\Qdet(d^* d)=\sum_B\Qdet B^*B$.
Here $B$ runs over choices of $n$ edges of $\G$.
If $B$ corresponds to a set of edges which has a component
with no cycle, that is, a component with $k-1$ edges and $k$
vertices, the corresponding matrix $B$ is singular (in fact has
a block form with nonsquare blocks) and so $\Qdet B=\sqrt{\det B}=0$.

Thus, as before, each component of $B$ must have exactly one cycle.
We now work directly with $B^*B$. Note that this is the
Laplacian on the subgraph defined by~$B$.
Suppose that $B$ has only one component; the general case is similar.
A~nonzero term in the expansion of the determinant corresponds
to a decomposition of $\sigma$ into cycles; the only such terms
which are nonzero are when adjacent vertices are paired,
each vertex is fixed, or the vertices on the unique cycle are
advanced or retreated in the direction of that cycle and the other vertices
are fixed.

Consider the case when some vertices are paired; if $i,j$ are paired in $B^*B$
then $B$ maps $i$ and $j$ to the same edge $e_{ij}$, so this term
does not contribute. If no terms are paired, the
expansion of the $Q$-determinant for $B^*B$ therefore has exactly
four terms.
It has two terms in which $\sigma$ fixes the bush vertices
and advances the cycle vertices around the cycle in one of the two orientations.
It also has two terms in which all vertices are fixed:
the vertices in the cycle advance to the edges in one of the two directions
and then retreat back to their original positions. The signs work out
as before.
\end{pf}
\begin{theorem}[(Noncommutative Cauchy--Binet formula)]\label{CB}
\[
\Qdet(M_1M_2) = \sum_{B_1} \Qdet(B_1B_2),
\]
where $B_1$ runs over maximal minors of $M_1$, and $B_2$ is the corresponding
minor of~$M_2$.
\end{theorem}
\begin{pf*}{Proof \textup{(Sketch)}} Suppose $M_1$ is $n\times m$ with $m>n$.
The standard proof only relies on the multilinearity of the determinant.
Write columns of $M_1M_2$
as linear combinations of columns of $M_1$
with coefficients in $M_2$ (acting on the right).
Use multilinearity to write this as a sum of $m^n$ $Q$-determinants;
the $Q$-determinant of one of these is zero if two of the same
columns are used.~Group the remaining $m$ choose $n$ nonzero terms
into products\break $\Qdet B_1\Qdet B_2$.
\end{pf*}

\section{Graphs on surfaces}

\subsection{Definitions}
Let $\Sigma$ be an oriented surface, possibly with boundary,
and $\G$ a graph embedded on $\Sigma$
in such a way that complementary components are contractible or
peripheral annuli (i.e., an annular neighborhood of a boundary
component).
We call the pair $(\G,\Sigma)$ a \textit{surface graph}.

For each boundary component there is a ``peripheral'' cycle
on $\G$, consisting of those vertices and edges bounding the same
complementary component as the boundary component.

An \textit{finite lamination} or \textit{simple closed curve system}
on $\Sigma$ is (the isotopy
class of) a finite set of pairwise disjoint simple closed curves on
$\Sigma$
each of which has nontrivial homotopy type (i.e., no curve bounds a
disk). We allow two or more curves in the lamination
to be isotopic.

\subsection{Bundles on surface graphs}

Given a surface graph $(\G,\Sigma)$ a vector bundle on $\G$
with connection $\Phi$ is
\textit{flat} if it has
trivial monodromy around any face (contractible complementary component
of $\Sigma$; peripheral cycles do not bound faces).
In this case, given a closed loop starting at a base point $v$,
the monodromy around a loop only depends
on the
homotopy class of the (pointed) loop in $\pi_1(\Sigma)$, and so the monodromy
determines a representation of $\pi_1(\Sigma)$
into $\operatorname{End}(V)$.
This representation depends on the base point for $\pi_1$;
choosing a different
base point will conjugate the representation.

Conversely, let $\rho\in\operatorname{Hom}(\pi_1(\Sigma
),\operatorname{End}(V))$
be a representation of $\pi_1(\Sigma)$ into $\operatorname{End}(V)$;
for a fixed base point there is a unique flat line bundle, up to gauge
equivalence,
with monodomy $\rho$. It is easy to construct: start with trivial
parallel transport on the edges of a spanning tree of $\G$. Now on
each additional edge, adding it to the tree makes a unique cycle;
define the parallel transport along this edge to be
the image of the homotopy class of this cycle under $\rho$. Thus, we
have a
correspondence between flat bundles modulo gauge equivalence and
a space $X=\operatorname{Hom}(\pi_1(\Sigma),\operatorname
{End}(V))/\operatorname{End}(V)$, the space of homomorphisms
of $\pi_1(\Sigma)$ into $\operatorname{End}(V)$ modulo conjugation.
This space $X$ has the structure of an algebraic variety: it can be represented
using variables for the matrix entries of a set of generators of $\pi
_1(\Sigma)$,
modulo an ideal corresponding to the relations in $\pi_1(\Sigma)$ and the
conjugation relations. The variety $X$ is called the
\textit{representation variety} of $\pi_1(\Sigma)$ in $\operatorname{End}(V)$.

When dealing with flat bundles it is natural to restrict our
attention to CRSFs each of whose cycles is topologically nontrivial,
since other CRSFs have zero weight in $\det\Delta$.
We call these \textit{incompressible CRSFs}.
The cycles in an incompressible CRSF form a finite lamination.

We now restrict to $V=\C^2$ and $\SU$-connections.
For a flat bundle, $\det\Delta$ is a regular function on the
representation variety
$X=\operatorname{Hom}(\pi_1(\Sigma),\SU)/\allowbreak\SU$
(i.e., it is a polynomial function of the matrix entries).
Remarkably, by varying the representation
it is possible to extract from $\det\Delta$ those terms for
incompressible CRSFs
whose cycles have any given set of homotopy types.
The following theorem is due to Fock and Goncharov \cite{FG}:
\begin{theorem}[(\cite{FG}, Theorem 12.3)]\label{theo11}
If $\Sigma$ has nonempty boundary,
the products
\[
\prod_{\mathrm{cycles}\ \gamma}\Tr(\gamma)
\]
over all finite laminations
form a basis for the vector space of regular functions on the
representation variety
$\operatorname{Hom}(\pi_1(\Sigma),\SU)/\SU$.
\end{theorem}

In our case,
$\det\Delta$ is a regular function on $X$;
hence it can be written
\[
\det\Delta=\sum_{L}c_L\prod_{\mathrm{cycles}}\operatorname{Tr}(w),
\]
where
$c_L\in\Z$, the sum is over finite laminations $L$,
the product is over curves $\gamma$ in $L$ and $w$ is the monodromy of
$\gamma$.
By an integral change of basis
(in fact an upper triangular integer matrix, if we order isotopy types of
finite laminations by inclusion)
we can write this as
\[
\det\Delta=\sum_L c'_L\prod_{\mathrm{cycles}}\bigl(2-\operatorname{Tr}(w)\bigr),
\]
and therefore by Theorem \ref{mainSL2}, $c'_L$ is the desired weighted
sum of CRSFs
with lamination type $L$.

The actual extraction of a coefficient can be done as follows.
There is~a~na\-tural measure $\mu$ on $X$, essentially just the product
of $k$ copies of Haar measure on $\SU$, where $k$ is the rank of the
free group $\pi_1(\Sigma)$.
The regular functions discussed above form a
basis for the Hilbert space $L^2(X,\mu)$ (not an orthonormal basis but
obtained from an orthonormal
basis by an upper-triangular linear transformation). Thus, one can
obtain $c_L$ above
by an integration of $\det\Delta$ against a particular function over~$X$.

\subsection{Examples}

\subsubsection{Annulus}
The simplest example of a surface graph, after a graph on a disk,
is a graph on an annulus.\vadjust{\eject}

Suppose $(\G,\Sigma)$ is a surface graph on an annulus.
Since $\pi_1(\Sigma)=\Z$ is Abelian, it will suffice to work with
a line bundle (a $\Z$ subgroup of $\SL$ will have the same traces
as a $\Z$ subgroup of diagonal matrices in $\SL$, which is equivalent
to taking a line bundle connection).
Let $z\in\C^*$ be the monodromy of a flat line bundle.

Then $P(z)=\det\Delta$ is a Laurent polynomial in $z$. It is
reciprocal: $P(z)=P(1/z)$ by Theorem \ref{mainline}.
\begin{theorem}\label{reality}
$P(z)$ is a reciprocal polynomial with real and positive roots
and a double root at $z=1$.
\end{theorem}

We emphasize here that the roots of $P$ are unrelated to the
eigenvalues of $\Delta$.
Rather they are special values of the monodromy $z$ for which $\Delta$
is singular.
We conjecture that the roots are distinct.
\begin{pf*}{Proof of Theorem \ref{reality}}
When $z=1$, the line bundle is trivializable and $\Delta$
is the usual graph Laplacian. Thus, $P(1)=0$. However, there exists a~CRSF on
$\G$ with a single component winding once around $\Sigma$: just\vspace*{1pt} take any
cycle of this type and complete it to a CRT. So in Theorem
\ref{mainline}, the coefficient
of $(2-z-1/z)^1$ is nonzero, and $z=1$ is a double root (i.e., not
of higher order).

We prove reality by a deformation
argument. At each real root $r>1$ of $P$, $\Delta$ has a kernel $W_r$.
We claim that $\dim W_r=1$.
If it were of larger dimension, let $f_1,f_2\in W_r$ be independent.
These are
real and harmonic
for $\Delta$, and so lift to actual harmonic functions on the embedded
graph $(\tilde\G,\tilde\Sigma)$ on the universal cover $\tilde
\Sigma$
of $\Sigma$. These functions have the property that
$f(x+1)=r f(x)$ where ``$+1$'' represents the deck transformation, that is,
on a path winding once around the annulus the values of $f_1,f_2$ are
multiplied by $r$.

Let $v$ be a vertex on the boundary of $\Sigma$ (i.e.,
on a peripheral cycle) and $f_3$ be a linear
combination of $f_1,f_2$ which is zero at $v$. Then $f_3$
is a harmonic function on $\tilde\G$ which is zero on a biinfinite
sequence of boundary points. We claim that this is impossible unless $f_3$
is identically zero. By the mean-value principle,
from each zero of a (not identically zero) harmonic function on an infinite
graph there are two paths to $\infty$, one on which the function is
increasing (and eventually positive) and another on which the function
is decreasing
(and eventually negative).
Let $\alpha$ be a choice of increasing and eventually positive
infinite path from a fixed lift $\tilde v$ of~$v$,
and $\beta$ be a choice of decreasing and eventually negative
infinite path
from the same point $\tilde v$.
The union of the lifts of $\alpha$ starting at all lifts of $v$ is a
nonnegative connected set of vertices
containing all lifts of $v$, and containing a strictly positive path
from $-\infty$ to $\infty$;
this contradicts the existence of $\beta$.

This completes the proof of the claim that $f_3$ must be identically zero,
and thus the claim that $\dim W_r=1$.

The roots of $P$
vary continuously with the conductances. We can find a set of
conductances for which all roots are real, see the next paragraph. As
we vary the conductances, if two real roots $r_1,r_2$ collide
and become a~complex conjugate pair,
the sum $W_{r_1}\oplus W_{r_2}$ at the point of
collision must be a two-dimensional subspace of the kernel of $\Delta$,
which we showed above cannot happen.
Thus, the roots must remain real.
From Theorem \ref{mainline}, $P(z_1)$ cannot have negative real
roots so all roots are real and positive.

We now find specific conductances for which all roots are real.
Take a~set of maximal cardinality of disjoint cycles,
each winding around $\Sigma$,
and complete
it to a CRSF. Let $C_1,\ldots,C_k$ be the cycles in order
from innermost to outermost. Add $k-1$ more edges to join these CRTs up into
a connected set~$U$ (a ``chain of loops,'' with bushes),
so that $C_i$ is connected to $C_{i+1}$ with a~path
$\gamma_i$. Put a
large conductance $R$ on every edge of $U$ except for
one edge of every $C_i$,
and one edge of every path $\gamma_i$;
these edges of $U$ have weight $1$.
Put a small conductance $\eps$ on the remaining edges of $\G$
not in $U$.

%
\begin{figure}

\includegraphics{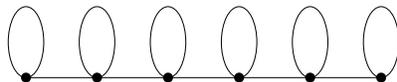}

\caption{A chain of loops. The roots of $\det\Delta$ are real and distinct.}\label{loopline}
\end{figure}

When $\eps$ is small and $R$ is large,
$\det\Delta$ is to leading order a power of $R$ times the determinant
of the graph of Figure \ref{loopline}. This is because with high
probability all
edges with conductance $R$, none of the edges
of conductance~$\eps$, and some of the edges of conductance $1$
will be present in a random CRSF. Contracting all edges of
conductance $R$ and removing those of conductance~$\eps$,
we are left with the graph of Figure \ref{loopline} with all edges of
conductance $1$.
For this graph, it is easy to verify that the roots $r_i$ of
$\det\Delta$ are real and distinct (see the example after Corollary
\ref{Bernoullicor}).
The actual Laplacian is a small perturbation of this
so its determinant also has real distinct roots.
%
%
%
%
%
\end{pf*}
\begin{cor}\label{Bernoullicor}
The number of cycles in a uniform random incompressible
CRSF on an annulus
is distributed as $1$ plus a sum of $k$ independent Bernoulli
random variables, where $k+1$ is the maximal number of
nonintersecting incompressible cycles one can simultaneously draw on
$\G$.
\end{cor}
\begin{pf}
Let $w=2-z-1/z$. Then $Q(w)=P(z)$ has only real negative roots:
\[
Q(w)=w\prod_{i=1}^k (w+\lambda_i).
\]
By Theorem \ref{mainline},
$Q(w)/Q(1)$ is the probability generating function for the number of loops.
It is also the probability generating function for a sum of independent
Bernoullis,
with the $i$th Bernoulli being biased as
$(\frac1{1+\lambda_i},\frac{\lambda_i}{1+\lambda_i})$.
\end{pf}

As an example,
let us compute, for a rectangular cylinder, $\det\Delta$ and the corresponding
distribution of cycles in a uniform random incompressible CRSF.
Let $\G_m$ be a line graph of length $m$ and $H_{m,n}=\G_m\times\Z_n$.
The eigenvalues of $\Delta_{\G_m}$ are $2+2\cos\frac{k\pi}{m}$ for
$k=1,2,\ldots,m$
[the corresponding eigenvectors are $f_k(x) = \cos\frac{\pi
k(x+1/2)}{m}$ for $x=0,1,2,\ldots,m-1$].
By (\ref{Ncrsfs}), we have
\begin{eqnarray*}
\det\Delta&=&(2-z-z^{-1})\prod\biggl(\Ch_n(\lambda+2)-z-\frac1z\biggr)
\\
&=&w\prod_{\lambda} w+\Ch_n(\lambda+2)-2,
\end{eqnarray*}
where $w=2-z-\frac1z$ and $\Ch_n$ is defined by
$\Ch_n(\alpha+\frac1{\alpha})=\alpha^n+\alpha^{-n}$ (a variant of
the Chebyshev polynomial).
For nonnegative $\lambda$ and large $n$, $\Ch_n(\lambda+2)$ is large
unless $\lambda$ is near $0$,
so the relevant roots are $\lambda_j=2+2\cos\frac{\pi(m-j)}m$ for
small $j$,
$j=1,2,\ldots.$
We have $2+\lambda_j=2+\frac{\pi^2j^2}{m^2}+O(\frac{j}m)^4 = \alpha
_j+1/\alpha_j$ where
$\alpha_j=1+\frac{\pi j}{m}+O(\frac{j}m)^2$.
Thus, in the limit $m,n\to\infty$ with $m/n\to\tau$ the roots satisfy
$\Ch_n(\lambda_j+2)=2\cosh\frac{\pi j}{\tau}+o(1)$.

The limit probability generating function for the number of cycles is then
\[
Q(w)/Q(1)=w\prod_{j=1}^\infty\biggl(\frac{w+2\cosh{(\pi
j}/{\tau})-2}{2\cosh({\pi j}/{\tau})-1}\biggr).
\]

For example, for a square annulus ($m=n$) the probability of a single
cycle is approximately
$95\%$.

\subsubsection{Torus example}

We compute the Laplacian determinant for
a flat line bundle on a $n\times n$ grid on a torus, that is, for the graph
$\Z^2/n\Z^2$. In this case, the determinant is
\[
\det\Delta=\prod_{\zeta^n=z_1}\prod_{\xi^n=z_2} 4-\zeta-\frac
1{\zeta}-\xi-\frac1{\xi}.
\]
This product was evaluated in \cite{BdT},
yielding
$\det\Delta=C_nP_n(z_1,z_2)$ where $C_n$ is a constant tending to $
\infty$ and $P_n$ tends to $P$,
where
%
\begin{equation}\label{P}
P(z_1,z_2) = \sum_{j,k\in\Z} e^{-({\pi
}/{2})(j^2+k^2)}z_1^jz_2^k(-1)^{j+k+jk}.
\end{equation}
We wish to rewrite this as
\[
P(z_1,z_2) = \sum_{(j,k)\ \mathrm{primitive}}\sum_{m\ge1} C_{jkm}
(2-z_1^jz_2^k-z_1^{-j}z_2^{-k})^m,
\]
where the first sum is over primitive vectors $(j,k)$ in $\Z^2$, one
per direction
[i.e., only one of $(j,k)$ and $(-j,-k)$ appears].
Thus, $C_{jkm}$ will be the relative probability that a CRSF has $m$ cycles
in homology class $(j,k)$.
To this end, pick a primitive vector $(j,k)$ and consider the terms in
(\ref{P})
with monomials $z_1^{mj}z_2^{mk}$ for $m\in\Z$. Letting
$u=z_1^jz_2^k$ and
$q=e^{-({\pi}/{2})(j^2+k^2)}$ the sum of these terms is
\[
\sum_{m\in\Z}q^{m^2}(-u)^m,
\]
which can be rewritten as a power series in $2-u-1/u$ as\footnote{We
use the fact that
\[
2-z^\ell-z^{-\ell}=\sum_{m\ge1}\frac{\ell}{m}\pmatrix{m+\ell
-1\cr2m-1}(-1)^{m+1}(2-z-z^{-1})^m.
\]
}
\[
\sum_{\ell\in\Z}q^{\ell^2}(-1)^\ell+ \sum_{m=1}^{\infty}\Biggl(
\sum_{\ell=1}^{\infty}\frac{\ell}{m}\pmatrix{m+\ell-1\cr
2m-1}(-1)^{\ell+m}q^{\ell^2}\Biggr)(2-u-1/u)^m .
\]
Thus, we find
\[
C_{jkm} =
\sum_{\ell=1}^{\infty}\frac{\ell}{m}\pmatrix{m+\ell-1\cr
2m-1}(-1)^{\ell+m}e^{-({\pi}/2)(j^2+k^2)\ell^2}.
\]

For example, the probability of a single cycle in homology class
$(1,0)$ is
\[
\frac{C_{101}}{\sum C_{jkm}}\approx41\%.\vspace*{15pt}
\]

\subsubsection{Pants example}

We consider here the case when $\Sigma$ is a $3$-holed sphere.
This is the simplest case where the noncommutativity of the connection
plays a role.
We consider a flat $\SL$-connection on $(\G,\Sigma)$.
Let $a,b$ and $c$ be the monodromies around the three
holes, which satisfy $abc=1$.
Since $\pi_1(\Sigma)$ is the free group on two
generators, there are no other relations. Note that a simple closed
curve on $\Sigma$
has only one of three possibly homotopy types: it must be homotopic to one
of the three boundary curves. A finite lamination up to homotopy is therefore
described by a triple $(i,j,k)\ne(0,0,0)$ of nonnegative integers,
where there are $i$ curves homotopic to $a$, $j$ homotopic to $b$, and
$k$ homotopic
to $c$. Theorem \ref{mainSL2} gives the following.
\begin{theorem} Let $X=2-\operatorname{Tr}(a), Y=2-\operatorname{Tr}(b), Z=2-\operatorname{Tr}(c)=2-\operatorname{Tr}(ab)$.
Then
$\det\Delta=P(X,Y,Z)=\sum_{\mathrm{CRSFs}} X^iY^jZ^k$
where $(i,j,k)$ is the number of cycles with homotopy type
$a,b,c$, respectively.
\end{theorem}

Note in particular that $P$ is a polynomial with nonnegative
coefficients. What polynomials occur? Is there an analog of Theorem
\ref{reality} in this setting?
We do not know. However, it is not difficult to show that the set of
exponents $(i,j,k)$
of the nonzero monomials of $P(X,Y,Z)$ form a convex set in $\Z^3$:
let~$\rho_{ab}$ be the length in terms of the number of vertices
of the shortest path from face $a$ (the face with monodromy $a$) to
face $b$
(the face with monodromy~$b$).
Similarly define $\rho_{ac}$ and~$\rho_{bc}$. Then we can
simultaneously embed
$i,j,k$ cycles of types $a,b,c$ if and only if
%
\begin{eqnarray}\label{3ineq}
i+j&\leq&\rho_{ab},\nonumber\\
i+k&\leq&\rho_{ac},\\
j+k&\leq&\rho_{bc}\nonumber.
\end{eqnarray}
This is easily proved by induction on $\rho_{ab},\rho_{ac},\rho_{bc}$.
So the Newton polyhedron of $P$ is defined by $i,j,k\ge0$, $(i,j,k)\ne
(0,0,0)$,
and (\ref{3ineq}).

\section{Planar loop-erased random walk}

In this section, we use the vector-bundle Laplacian to
study the loop-erased walk on a planar graph.
In Section~\ref{onehole}, we compute the probability that the LERW
between two boundary points passes left or right of a given face.
In Section \ref{twoholes}, we compute this probability for two faces.

\subsection{LERW around a hole}\label{onehole}

Let $\G$ be a finite planar graph.
Let $z_1,z_2$ vertices on the outer boundary and $f$ a
bounded face.
We compute the probability that the LERW from $z_1$ to $z_2$
goes left of $f$. This can be computed quite simply using duality: given
a UST on a planar graph $\G$, the duals of the edges not in $T$ form a
UST of the planar dual graph. In the dual graph, our question is
equivalent to asking whether the
LERW from the face $f$ first reaches the outer face along the edge
between $z_1$
and $z_2$ or outside this edge. This is just the harmonic measure of
the path
from $z_1$ to $z_2$ as seen from $f$.

Here we will compute this in another way, using line bundle technology.
It illustrates
a general method which will be useful in the next section.

\begin{figure}

\includegraphics{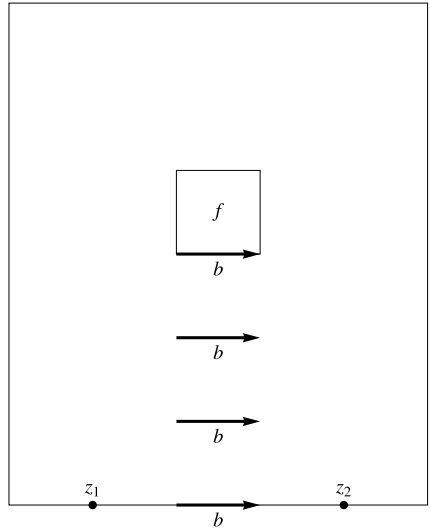}

\caption{Defining the connection for the LERW around a hole.}\label{zipper}
\end{figure}

Take a line bundle with monodromy $b$ around $f$
and trivial around all other faces.
This can be achieved by
putting in a zipper of edges (dual of a~simple path in the dual graph
from $f$ to the boundary) each with parallel transport $b$
as in Figure \ref{zipper}.
Suppose without loss of generality
that the zipper ends between $z_1$ and $z_2$ as in
the figure.

%

Put in an extra edge from $z_1$ to $z_2$ with parallel transport $\phi
_{z_1z_2}=a$. Then by Theorem \ref{mainline},
%
\begin{eqnarray}\label{p1p2p3}
\det\Delta&=&p_1 + p_2\biggl(2-a-\frac1a\biggr)+p_3
\biggl(2-\frac{a}{b}-\frac{b}{a}\biggr)
\nonumber\\[-8pt]\\[-8pt]
&=&p_1+2p_2+2p_3-a\biggl(p_2+\frac{p_3}b\biggr)-\frac1a(p_2+bp_3).\nonumber
\end{eqnarray}
Here $p_2=p_2(b)$
is the weighted sum of CRSFs which contain a cycle
containing
edge $a$ and going left of $f$ (when oriented so that $a$ is the edge
before~$z_1$); $p_3=p_3(b)$ is the weighted sum of CRSFs containing
a cycle through $a$ going ``right'' of~$f$.

We wish to compute $p_2(b=1), p_3(b=1)$.
These can be extracted from the coefficients of $a$ and $\frac1a$
in $\det\Delta$.

Letting $z_1,z_2$ be the first two vertices, we have
\[
\Delta= \Delta(b) + \pmatrix{1&-a&\bz\cr-a^{-1}&1&\bz\cr
\bz&\bz&\bz},
\]
where $\Delta(b)$ is the Laplacian without the extra edge from $z_1$
to $z_2$
(and $\bz$ represents a vector or matrix of zeros).
As a consequence,
\[
\det\Delta= \det\Delta(b)\det\left(I+\pmatrix{1&-a\cr-1/a&1
}\pmatrix{G(z_1,z_1)&G(z_1,z_2)\cr
G(z_2,z_1)&G(z_2,z_2)
}\right),
\]
where $G=\Delta(b)^{-1}$.

This last determinant is $1+G(z_1,z_1)+G(z_2,z_2)-a G(z_2,z_1)-\frac
{G(z_1,z_2)}{a}$.
Comparing coefficients with
(\ref{p1p2p3}) yields
\[
p_2=\frac{b^2G(z_2,z_1)-G(z_1,z_2)}{b^2-1}\det\Delta(b)
\]
and
\[
p_3=\frac{b(G(z_1,z_2)-G(z_2,z_1))}{b^2-1}\det\Delta(b).
\]

The desired quantity is
%
\begin{equation}\label{oneholeprob}
\frac{p_2}{p_2+p_3}=\frac{b^2
G(z_2,z_1)-G(z_1,z_2)}{(b-1)(bG(z_2,z_1)+G(z_1,z_2))}
\end{equation}
in the limit $b\to1$.

We need to compute $G(v,v')$ when $b$ is near $1$.
\begin{lemma}\label{oneterm} When $b=1+\eps$,
$G(v,v')=\frac{\kappa}{\det\Delta(b)}(1+\eps X(v,v')+O(\eps^2))$
where $\kappa$ is the number of spanning trees of $\G$ and
$X(v,v')$ is the expected signed number of crossings of the zipper
on a simple random walk from $v$ to~$v'$.
\end{lemma}
\begin{pf}
For a matrix $M$, let $M_A^B$ denote $(-1)^{i_A+i_B}$ times
the determinant of $M$ upon
removing rows $A$ and columns $B$.
Here $i_A$ is the sum of the indices of elements of $A$ and likewise
for $i_B$. (We can assume without loss of generality
that all relevant vertices have even index
so that we can ignore the signs.) Recall that for the standard
Laplacian $\Delta_0$, we have
for any $v,v'$ that $(\Delta_0)_v^{v'}= \kappa$, the number of
spanning trees of $\G$.
We also have $(\Delta_0)_{ab}^{cd}=G(a,c)-G(a,d)-G(b,c)+G(b,d)$\vspace*{2pt} where
$G$ is the
Green's function.\footnote{More generally,
\[
(\Delta_0)_{a_1,\ldots,a_k}^{b_1,\ldots,b_k}=[t]\kappa\det
\bigl(t+G(a_i,b_j)\bigr)_{1\le i,j\le k},
\]
where $[t]P(t)$ denotes the coefficient of $t^1$ in $P(t)$. This can be
proved by induction on $k$
and expanding the determinant on the left-hand side along row $a_k$.}

Letting $Z=\Delta(b)-\Delta_0$,
a matrix supported on the zipper edges.
We have
%
\begin{eqnarray}\label{detexpansion}
&&G(v,v')\det\Delta(b)\nonumber\\
&&\qquad= \Delta(b)_{v'}^{v}
=(\Delta_0+Z)_{v'}^{v}
=(\Delta_0)_{v'}^{v}+\sum_{x,y}Z_{x}^{y}(\Delta
_0)_{v',x}^{v,y}+\cdots
\nonumber\\
&&\qquad=\kappa+\sum_{\vec{uw}}(b-1)(\Delta_0)_{v',u}^{v,w}+\biggl(\frac
1b-1\biggr)(\Delta_0)_{v',w}^{v,u}+O(\eps^2)\\
&&\qquad=
\kappa\biggl(1+\eps\sum_{\vec
{uw}}G_0(v',u)+G_0(v,w)\nonumber\\
&&\qquad\quad\hspace*{12pt}{}-G_0(v',w)-G_0(v,u)+O(\eps^2)\biggr).\nonumber
\end{eqnarray}
Here the sum is over the zipper edges.
For an edge $\vec{uw}$ the quantity $G_0(v',u)+G_0(v,w)-G_0(v',w)-G_0(v,u)$
is the transfer current, equivalently, the
expected signed number of crossings of $uw$ of the
simple random walk started at $v'$ and stopped at $v$.
\end{pf}

By the lemma, when $b=1+\eps$ the probability (\ref{oneholeprob})
becomes
\[
\frac{p_2}{p_2+p_3}=1+\frac{X(z_2,z_1)-X(z_1,z_2)}2+O(\eps
)=1-X(z_1,z_2)+O(\eps)
\]
since $X(z_1,z_2)=-X(z_2,z_1)$.
\begin{cor}
The probability that the LERW from $z_1$ to $z_2$
goes left of $f$ is
the total amount of current flowing left of $f$ when $1$ unit of
current enters at $z_1$ and exits at $z_2$.
\end{cor}

\subsection{LERW around two holes}\label{twoholes}

Let $f_1,f_2$ be two faces in a finite planar graph $\G$ and
$z_1,z_2$ two points on the boundary.
Consider the LERW from $z_1$ to $z_2$.
We compute the probability that it goes left of both faces.
Although this can be computed in principle using Theorem \ref{theo11}
and the comments thereafter, in practice that method involves an integration
of $\det\Delta$ over the representation variety $X$; in the case of a
large graph
one might not have access to an explicit expression for $\det\Delta$.
Here, we
do the computation using only a bundle which is close to the trivial bundle;
the determinant as a~function on $X$
can be expanded as a power series around the trivial bundle
case and we only need the first two terms in the expansion, which we can
write down explicitly in terms of the standard green's function.

\begin{figure}[b]

\includegraphics{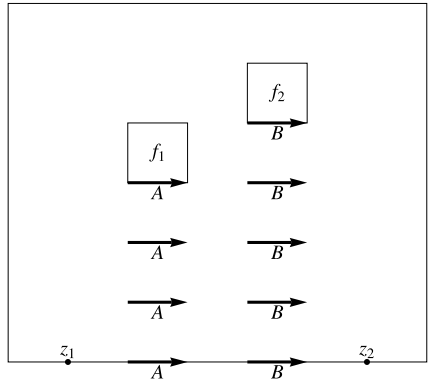}

\caption{Defining the connection for the LERW around two holes.}\label{twozips}
\end{figure}

We use an $\SL$-bundle.
Choose a bundle with mondromy $A$ around $f_1$ and $B$ around $f_2$,
and trivial monodromy around the other faces. This can be achieved by putting
trivial parallel transport on all edges except for two zippers,
one from each of $f_1,f_2$ to the boundary, as in Figure \ref{twozips}.
Let $Z_{A},Z_B$
denote these two zippers.

Add a new edge connecting $z_1$ to $z_2$ with parallel transport $C$.
Suppose that $A,B$ are close to the identity. Then a CRSF
will have with high probability only a loop
containing $C$. The first-order correction to this consists
in CRSFs having two loops, one of them containing $C$ and the other
of monodromy $A,B$ or $AB$ (or their inverses).
The possible loops in a CRSF with one or two loops,
one of which contains $C$, have
one of the seven monodromy
types (or their inverses):
\[
A,B,AB,C,CB^{-1}A^{-1},C(BA)^nA^{-1}(BA)^{-n}, C(AB)^nB^{-1}(AB)^{-n},
\]
where $n\in\Z$.

Let
\begin{eqnarray*}
A&=&\pmatrix{1+xu&y\sqrt{u}\cr z\sqrt{u}&1+wu},\\
B&=&\pmatrix{1+au&b\sqrt{u}\cr c\sqrt{u}&1+du},\\
C&=&\pmatrix{e&1\cr-1&0}
\end{eqnarray*}
be the representations in $\SL$, where $u$ is small and $x,y,z,w,a,b,c,d,e$
are variables. The equations
$\det A=\det B=1$ yield
$x+w-yz=o(u)$ and $a+d-bc=o(u)$.

Computing the traces of the above seven classes, we have
%
\begin{eqnarray}\label{7traces}
2-\Tr(A)&=&-(x+w)u+o(u),\nonumber\\
2-\Tr(B)&=&-(a+d)u+o(u),\nonumber\\
2-\Tr(AB)&=&-(a+d+w+x+cy+bz)u+o(u),\nonumber\\
2-\Tr(C)&=&2-e,\\
\qquad2-\Tr(C(AB)^nB^{-1}(AB)^{-n})&=&2-e\bigl(1+du+o(u)\bigr),\nonumber\\
2-\Tr(C(BA)^nA^{-1}(BA)^{-n})&=&2-e\bigl(1+wu+o(u)\bigr),\nonumber\\
2-\Tr(CB^{-1}A^{-1})&=&2-e\bigl(1+(d+w+bz)u\bigr).\nonumber
\end{eqnarray}

The relevant quantities for us are the last four and the
products of any one the first three of these with the last four.

We let $N_{ll}$ be the weight sum of configurations in which the LERW
from $z_1$ to $z_2$ goes left of both faces, and having no other loops
(i.e., CRSFs with one loop containing the extra edge, with the property
that this loop surrounds both $f_1,f_2$).
Similarly, define $N_{lr},N_{rl},N_{rr}$. Let $N_1,N_2,N_3$ be,
respectively, the weighted sum of those configurations
having one loop containing $C$ and one loop surrounding,
respectively, hole $f_1$, hole $f_2$,
and both holes.

If we extract the coefficient of $e$ in $-\Qdet\Delta$, from (\ref
{7traces}) the
result will be of the form
\begin{eqnarray*}
&&N_{ll}+N_{lr}(1+du)+N_{rl}(1+wu)+N_{rr}\bigl(1+(d+w+bz)u\bigr)
\\
&&\qquad{} -N_1(w+x)u-N_2(a+d)u-N_3(a+d+w+x+cy+bz)u+O(u^2).
\end{eqnarray*}

Recall that these variables are subject to the constraints
$x+w-yz=o(u)$ and $a+d-bc=o(u)$. The coefficients of the $N_*$
are linearly independent given these constraints.
Thus given an explicit expression for
$\Qdet\Delta$ to first order in $u$, it is
a simple matter
to extract these various coefficients $N_*$. The desired probabilities
are
\[
p_{ll}=\frac{N_{ll}}{N_{ll}+N_{lr}+N_{rl}+N_{rr}}
\]
and so on.

Let us show how to compute the coefficient of $e$
in $\Qdet\Delta$, to first order in $u$.

Indexing the vertices so that $z_1,z_2$ are the first two vertices, we have
\[
\Delta(A,B,C)= \Delta(A,B) + \pmatrix{I&-C&\bz\cr-C^{-1}&I&\bz\cr
\bz&\bz&\bz},
\]
where $\Delta(A,B)$ is the Laplacian
without the extra edge from $z_1$
to $z_2$
(and~$\bz$ represents a vector or matrix of zeros).
As a consequence,
\begin{eqnarray*}
&&\det\Delta(A,B,C) \\
&&\qquad= \det\Delta(A,B)\det\left(I+\pmatrix
{1&0&-e&-1\cr0&1&1&0\cr0&1&1&0\cr
-1&-e&0&1}\right.\\
&&\hspace*{129pt}{}\times\left.\pmatrix{G\bigl(z_1^{(1)},z_1^{(1)}\bigr)&\cdots
&&G\bigl(z_1^{(1)},z_2^{(2)}\bigr)\cr
G\bigl(z_2^{(2)},z_1^{(1)}\bigr)&\cdots&&G\bigl(z_2^{(2)},z_2^{(2)}\bigr)
}\vphantom{I+\pmatrix
{1&0&-e&-1\cr0&1&1&0\cr0&1&1&0\cr
-1&-e&0&1}}\right),
\end{eqnarray*}
where $G=\Delta(A,B)^{-1}$.
Note that $G$ is self-dual,
since it is the inverse of a~self-dual matrix.
Thus, we have
\begin{eqnarray*}
\Qdet\Delta(A,B,C)
&=&\Qdet\Delta(A,B)\\
&&{}\times \Qdet \left(I+\pmatrix{
I&-C\cr-C^{-1}&I}\pmatrix{
\mathbf{G}(z_1,z_1)&\mathbf{G}(z_1,z_2)\cr\tilde\mathbf{G}(z_1,z_2)&
\mathbf{G}(z_2,z_2)}\right),
\end{eqnarray*}
where $\mathbf{G}(z_1,z_1)$ and $\mathbf{G}(z_2,z_2)$ are diagonal
matrices.
The right-hand side is $\Qdet\Delta(A,B)$ times
\[
1+G\bigl(z_1^{(1)},z_1^{(1)}\bigr)-G\bigl(z_1^{(1)},z_2^{(2)}\bigr)+G\bigl(z_1^{(2)},z_2^{(1)}\bigr)+G\bigl(z_2^{(1)},z_2^{(1)}\bigr)-e
G\bigl(z_1^{(2)},z_2^{(2)}\bigr).
\]

Thus, the coefficient of $e$ in $-\Qdet\Delta(A,B,C)$ is
\[
G\bigl(z_1^{(2)},z_2^{(2)}\bigr)\Qdet\Delta(A,B)
\]
to first order in $u$.\vadjust{\eject}

It is possible to give an exact expression for
$G(z_1^{(2)},z_2^{(2)})\Qdet\Delta(A,B)$
as a sum over $Z_A$ and $Z_B$ of products of standard Green's functions.
Since $A,B$ are close to the identity,
we can write
\[
\Delta(A,B) = \pmatrix{\Delta_0&0\cr0&\Delta_0}
+ Z,
\]
where $Z$ is a perturbation supported on the zippers.
As a consequence, both
\[
G\bigl(z_1^{(2)},z_2^{(2)}\bigr)\det\Delta(A,B)=\det\Delta(A,B)_{z_2^{(2)}z_1^{(2)}}
\]
and
$\Qdet\Delta(A,B)$ itself have expansions in powers of $u$. The ratio
of these
is the desired quantity.

For the present purposes, we need the expansion of
$G(z_1^{(2)},z_2^{(2)})\det\Delta(A,B)$ to first order (which is the
term of order $u^2$) and $\sqrt{\det\Delta(A,B)}$
to first order (which is the term of order $u$).
Unfortunately these computations, although elementary, can be quite long,
and we have not been able to get an explicit answer for a main
case of interest which is the upper half plane.\vspace*{12pt}

\section{Green's function and random walks}

We give here a probabilistic interpretation of the Green's function of
the line bundle
Laplacian in the case where each edge has parallel transport either $1$
or $z$.
Let $\G$ be a graph,~$E_Z$ a collection of directed edges, and consider
the line bundle with parallel transport $z$ on $E_Z$ and $1$ elsewhere.
Let $G$ be the inverse of the line bundle Laplacian.

Let $v,v'$ be vertices of $\G$.
Let us first ask: what is the distribution
of the signed number of crossings of $E_Z$ on the SRW from $v$ to $v'$?

Let $P$
be the transition matrix of the simple random walk on $\G$ (ignoring
the line bundle)
with absorbing state $v'$. Thus, $P_{v_2v_1}$ is the probability of
going to $v_2$
from state $v_1$ (as usual the indices are reversed so that composition
respects the natural path order).
Multiply the entry $P_{v_2v_1}$ by $z$ or $1/z$ if $v_1v_2$ is in $E_Z$
or $-E_Z$,
respectively.
With $v'$ as the vertex with
last index~$P$ has the form $P=\bigl({{N\atop A} \enskip{0\atop1}}\bigr)$.
Then the probability generating function of the number of crossings
is
\[
\sum_{k=0}^\infty(AN^k)_{v'v}=\bigl(A(I-N)^{-1}\bigr)_{v'v}.
\]
\begin{prop}\label{Gvv'}
The probability generating function of the signed
number of crossings of $E_Z$ on a simple random walk
from $v$ to $v'$ is
\[
\bigl(A(I-N)^{-1}\bigr)_{v'v}=\frac{G(v,v')}{G(v',v')}.
\]
\end{prop}
\begin{pf}
Note that $((I-P)D)_{v_1v_2}=\Delta_{v_1v_2}$ for $v_2\ne v'$,
where $D$ is the diagonal matrix of vertex degrees.
As a consequence
\begin{eqnarray*}
\bigl(A(I-N)^{-1}\bigr)_{v'v}&=&\sum_{v_2}A_{v'v_2}(I-N)^{-1}_{v_2v} \\
&=&\sum_{v_2}A_{v'v_2}\frac{\cof(I-N)_{v}^{v_2}}{\det(I-N)}\\
&=&\sum_{v_2}A_{v'v_2}\frac{\cof(\Delta D^{-1})_{vv'}^{v_2v'}}{
\cof(\Delta D^{-1})_{v'}^{v'}}\\
&=&\sum_{v_2}A_{v'v_2}\frac{\cof\Delta_{vv'}^{v_2v'}(
{\deg(v')\deg(v_2)}/{\det D})}{
(\det\Delta)G(v',v'){\deg(v')}/{\det D}}\\
&=&\frac1{(\det\Delta)G(v',v')}\sum_{v_2}\Delta_{v'v_2}\cof\Delta
_{vv'}^{v_2v'}\\
&=&\frac{G(v,v')}{G(v',v')}.
\end{eqnarray*}
\upqed
\vspace*{6pt}
\end{pf}

The quantity $1-\frac1{G(v',v')\deg v'}$ also has a probabilistic
interpretation.
Choosing a gauge such that $\phi_{v'v}=1$ for all $v$ adjacent to $v'$,
and using
\[
\sum_{v_2}\Delta_{v,v_2}G(v_2,v') = \delta_{v,v'}
\]
we have (setting $v=v'$)
\[
\deg v' G(v',v')-\sum_{v_2\sim v'}G(v_2,v')=1
\]
giving
\[
1-\frac1{G(v',v')\deg v'}=\sum_{v_2\sim v'}\frac{1}{\deg v'}\frac
{G(v_2,v')}{G(v',v')}.
\]
Thus, we have the following proposition.
\begin{prop}\label{Gvv}
$1-\frac1{G(v',v')\deg v'}$ is the probability generating function
of the signed number of crossings
of $E_Z$ on a SRW started at $v'$ and stopped on its first return
to $v'$.
\end{prop}

Now let $z=1+\eps$ with small $\eps$ and
%
\begin{equation}\label{Gvv'expansion}
\frac{G(v,v')}{G(v',v')}=1+X_1\eps+X_2\eps^2+\cdots,
\end{equation}
so that $X_1$ is the expected signed number of crossings
of $E_Z$ on a SRW from~$v$ to~$v'$. Similarly let
%
\begin{equation}\label{Gv'v'expansion}
1-\frac1{G(v',v')\deg v'}=1+Y_2\eps^2+Y_3\eps^3+\cdots,
\end{equation}
so that $2Y_2$ is the expectation of the square of the
signed number
of crossings of $E_Z$ of the SRW from $v'$ until its first return to $v'$
(note that $Y_1$ is zero and $Y_3=-Y_2$ by symmetry).
\begin{prop} The determinant of the line bundle Laplacian $\Delta$
above satisfies
%
\begin{equation}\label{detdelta1}
\det\Delta= \kappa\deg v'\bigl(-Y_2\eps^2+Y_2\eps^3 + O(\eps^4)\bigr),
\end{equation}
where $2Y_2$ is the expected square of the signed number of crossings
of $E_Z$ of a SRW from $v'$ until its first return to $v'$.
In particular $Y_2\deg v'$ does not depend on~$v'$.
\end{prop}
\begin{pf}
From Lemma \ref{oneterm}, we have
\[
G(v,v')\det\Delta= \kappa\bigl(1+X_1\eps+O(\eps^2)\bigr).
\]
However by (\ref{Gvv'expansion}) and
(\ref{Gv'v'expansion}), we have
\[
G(v,v')\det\Delta=\frac1{\deg v'}\frac{1+X_1\eps+X_2\eps^2+O(\eps
^3)}{-Y_2\eps^2+Y_2\eps^3+O(\eps^4)}\det\Delta,
\]
so we must have $\det\Delta=\kappa\deg v'(-Y_2\eps^2+Y_2\eps
^3+O(\eps^4))$.
\end{pf}

We can get a full expansion of
$\det\Delta$ as a power series in $z$ around $z=1$ as follows.
We have
\[
\frac{\partial}{\partial z}\log\det\Delta=
\sum_{uw\in E_Z}G_{u}^{w}-\frac1{z^2}G_{w}^{u}.
\]
This follows from differentiating the entries in $\Delta$.
Here using Propositions~\ref{Gvv} and~\ref{Gvv'},
the quantities on the right-hand side have
explicit probabilistic interpretations.
Integrating both sides from $z=1$ to $z=z_0$, gives the desired
expansion.

%
%
%
%
%

\section{Monotone lattice paths}\label{monotone}

Let $\Gamma_{m,n}$ be the $m\times n$ grid
on a torus [the nearest neighbor grid $\{0,1,\ldots,m\}\times\{
0,1,\ldots,n\}$
with opposite sides identified: $(m,j)\sim(0,j)$ and $(i,n)\sim(i,0)$].
Choose for each vertex a north- or east-going edge, independently and with
probability $p=1/2$. The resulting configuration of edges makes a
directed CRSF, and we wish to determine the distribution of the number and
homology type of the cycles. This problem was studied in \cite{HK} who showed
among other things that
when $m/n$ is close to a rational with small denominator the number of
cycles is (in the limit $m,n\to\infty$) tending to a Gaussian with
expectation on the order
of $\sqrt{n}$.

Here we show how to compute explicitly for each $m,n$
the probability generating function of the total homology class of the cycles
(which determines both the number of cycles and their direction).

We make $\Gamma_{m,n}$ into a directed graph, directing all edges
northward or eastward.
Let $\Phi$ be a flat line bundle on $\Gamma_{m,n}$ with
monodromy $z$ and $w$ for loops in homology class $(1,0)$ and $(0,1)$,
respectively.
$\Phi$ can be obtained by putting parallel transports $1$ on all edges
except edges $(m-1,i)(0,i)$ which get parallel transport $z$ and edges
$(i,n-1)(i,0)$ which get parallel transport~$w$.

By Theorem \ref{directed}, the determinant of the line bundle Laplacian
is
%
\begin{equation}\label{Fzw}
F_{m,n}(z,w)=\det\Delta= \sum_{j,k\ge0} C_{j,k}(1-z^pw^q)^\ell,
\end{equation}
where $(j,k)=(p\ell,q\ell)$ and $p,q$ are coprime, and $C_{j,k}$ is the
number of CRSFs with total homology class $(j,k)$, that is, with $\ell
$ cycles
of homology $(p,q)$.

This determinant can be explicitly computed using a standard Fourier
diagonalization of $\Delta$.\vspace*{6pt}
\begin{prop}
%
\begin{equation}\label{Fprod}
F_{m,n}(z,w) = \prod_{u^m=z}\prod_{v^n=w}2-u-v.\vspace*{3pt}
\end{equation}
\end{prop}
\begin{pf}
If $u^m=z$ and $v^n=w$, then the function $f(x,y)=u^xv^y$ is an
eigenvector of $\Delta$
with eigenvalue $2-u-v$. These eigenvectors are independent and span
$\C^{mn}$.\vspace*{3pt}
\end{pf}

From (\ref{Fprod}), we can extract the coefficients $C_{m,n}$.
Let us consider for simplicity the case $m=n$. In \cite{HK}, it was
shown that
with probability tending to~$1$ as $n\to\infty$,
all cycles will have homology class $(1,1)$.
Thus, from (\ref{Fzw}),
%
\begin{equation}\label{2sums}
F_{n,n}(z,w) = \sum_{j\ge0}C_{j,j}(1-zw)^j + \sum_{j\ne k\ge
0}C_{j,k}(1-z^pw^q)^\ell,
\end{equation}
where the second sum is negligible, in the sense that the sum of its
coefficients
is $o(\sum C_{j,j})$. We can thus ignore the second sum, let $z=1$ and expand
$F_{n,n}(1,w)$ around $w=1$.
We have
\[
F_{n,n}(z,w) = \prod_{u^n=z,v^n=w} 2-u-v=\prod_{u^n=z}(2-u)^n-w.
\]
Letting $Y=1-w$ and $F_{n,n}(1,w)=H_n(1,Y)$, we have
\[
H_n(1,Y)=\prod_{j=0}^{n-1}(2-e^{2\pi ij/n})^n-(1-Y)
\]
and
\[
\frac{H_n(1,Y)}{H_n(1,1)}=\prod_{j=-\lfloor n/2\rfloor+1}^{\lfloor
n/2\rfloor}
\frac{(2-e^{2\pi ij/n})^n-(1-Y)}{(2-e^{2\pi ij/n})^n}.
\]
The product $(2-e^{2\pi ij/n})^n$ is large in absolute value
unless $|j|$ is small, that is, when $|j|$ is
$O(n^{1/2})$, so the dependence
on $Y$ comes from terms with $|j|$ small. That is, the product
can be written as
\[
\prod_{|j|<n^{1/2+\eps}}(\cdots)\prod_{|j|>n^{1/2+\eps}}(\cdots)
\]
and the second product is $1+o(e^{-n^{\eps}})$.
When $|\theta|$ is small,
\[
(2-e^{i\theta})^n=(1-i\theta+\theta^2/2+\cdots)^n = e^{-in\theta
+n\theta^2+\cdots}.
\]
Plugging in $\theta=2\pi j/n$ yields
\[
\frac{H_n(1,Y)}{H_n(1,1)}\bigl(1+o(e^{-n^\eps})\bigr)=\prod_{j=-\infty
}^{\infty}\frac{e^{(2\pi j)^2/n}-1+Y}{e^{(2\pi j)^2/n}},
\]
where we have extended the range of $j$ without loss of precision.
Thus up to exponentially small errors,
\[
\frac{H_n(1,Y)}{H_n(1,1)}=Y\prod_{j\ne0} \bigl(1-e^{-(2\pi
j)^2/n}+Ye^{-(2\pi j)^2/n}\bigr)
\]
is the probability generating function for the number $N_{(1,1)}$ of
$(1,1)$-cycles.
We see that the number of $(1,1)$ cycles is a sum of independent
Bernoulli random
variables with the $j$th having bias $p_j=e^{-(2\pi j)^2/n}$.
The expectation of $N_{(1,1)}$ is then (using $t=2\pi j/\sqrt{n}$)
\[
\sum p_j = \sum_j e^{-(2\pi j)^2/n} \approx
\frac{\sqrt{n}}{2\pi}\int_{-\infty}^{\infty}e^{-t^2}\,dt=\sqrt
{\frac{n}{4\pi}}.
\]
The variance is
\[
\sum p_j(1-p_j)\approx\frac{\sqrt{n}}{2\pi}\int_{-\infty}^{\infty
}e^{-t^2}(1-e^{-t^2})\, dt=\sqrt{\frac{n}{4\pi}}\biggl(1-\frac1{\sqrt{2}}\biggr).
\]
Since the variance tends to $\infty$, the central
limit theorem implies that distribution
tends to a Gaussian as $n\to\infty$.

A similar computation holds when $m=pk, n=qk$ and $k\to\infty$.

\section{Open questions}
\begin{enumerate}
\item
Consider a graph embedded on a Riemann surface, such that edges are
geodesic segments. Now
use a naturally associated connection,
for example, the Levi--Civita connection on the tangent bundle.
What are the probabilistic consequences of choosing such a connection?

\item
Is there any combinatorial meaning to the vector bundle Laplacian
for higher-rank bundles?

\item What information does the line bundle Laplacian give about the
three-dimen\-sional
monotone lattice path problem, generalizing the results of Section \ref
{monotone}?
\end{enumerate}

\section*{Acknowledgments}
We thank Alexander Goncharov, Adrien Kassel and David Wilson
for important contributions and the referee
for helpful feedback.


%
\printaddresses

\end{document}